\documentclass[11pt,a4paper,reqno]{amsart} 
\usepackage{amsfonts,amsthm, amscd, epsfig, amsmath, amssymb,enumerate}

\numberwithin{equation}{section}
\DeclareMathSymbol{\leqslant}{\mathalpha}{AMSa}{"36} 
\DeclareMathSymbol{\geqslant}{\mathalpha}{AMSa}{"3E} 
\DeclareMathSymbol{\eset}{\mathalpha}{AMSb}{"3F}     
\renewcommand{\leq}{\;\leqslant\;}                   
\renewcommand{\geq}{\;\geqslant\;}                   
\newcommand{\inftwo}[2]{\inf_{\substack{#1 \\ #2}}} 

\newcommand{\grad}{\nabla} 

\newcommand{\la}{\label} 
 
\newcommand{\be}{\begin{equation}}

\def\1{\ifmmode {1\hskip -3pt \rm{I}} \else {\hbox {$1\hskip -3pt \rm{I}$}}\fi}

\newtheorem{Th}{Theorem}[section] 
\newtheorem{Le}[Th]{Lemma}

\newcommand{\cF}{\ensuremath{\mathcal F}} 
 
\newcommand{\cH}{\ensuremath{\mathcal H}}

\newcommand{\cK}{\ensuremath{\mathcal K}} 
\newcommand{\cL}{\ensuremath{\mathcal L}}

\newcommand{\cT}{\ensuremath{\mathcal T}}

\newcommand{\bbR}{{\ensuremath{\mathbb R}} }

\newcommand{\bbZ}{{\ensuremath{\mathbb Z}} }

\newcommand{\si}{\sigma} 
 
\newcommand{\var}{{\rm Var} }

\let\a=\alpha \let\b=\beta   \let\d=\delta  \let\e=\varepsilon
 \let\g=\gamma       \let\l=\lambda
      \let\o=\omega      
\let\r=\rho      
  
   \let\G=\Gamma  \let\L=\Lambda 
\let\O=\Omega      

\def\\{\hfill\break}

\def\thsp{\thinspace}

\def\tthsp{\kern .083333 em}

\def\?{\mskip -10mu}
\def\indbox#1{\hbox to \parindent{\hfil\ #1\hfil} }

\newfam\msafam
\newfam\msbfam
\newfam\eufmfam
\def\hexnumber#1{%
  \ifcase#1 0\or 1\or 2\or 3\or 4\or 5\or 6\or 7\or 8\or
  9\or A\or B\or C\or D\or E\or F\fi}
\font\tenmsa=msam10
\font\sevenmsa=msam7
\font\fivemsa=msam5
\textfont\msafam=\tenmsa
\scriptfont\msafam=\sevenmsa
\scriptscriptfont\msafam=\fivemsa        
\edef\msafamhexnumber{\hexnumber\msafam}%
\mathchardef\restriction"1\msafamhexnumber16
\mathchardef\ssim"0218
\mathchardef\square"0\msafamhexnumber03
\mathchardef\eqd"3\msafamhexnumber2C
\def\QED{\ifhmode\unskip\nobreak\fi\quad
  \ifmmode\square\else$\square$\fi}            
\font\tenmsb=msbm10
\font\sevenmsb=msbm7
\font\fivemsb=msbm5
\textfont\msbfam=\tenmsb
\scriptfont\msbfam=\sevenmsb
\scriptscriptfont\msbfam=\fivemsb

\font\teneufm=eufm10
\font\seveneufm=eufm7
\font\fiveeufm=eufm5
\textfont\eufmfam=\teneufm
\scriptfont\eufmfam=\seveneufm
\scriptscriptfont\eufmfam=\fiveeufm

\def\({\left(}
\def\){\right)}
\let\neper=e
\let\ii=i

\def\nep#1{ \neper^{#1}}

\def\tc{\thsp | \thsp}
\let\<=\langle
\let\>=\rangle

\outer\def\nproclaim#1 [#2]#3. #4\par{\medbreak \noindent
   \talato(#2){\bf #1 \Thm[#2]#3.\enspace }%
   {\sl #4\par }\ifdim \lastskip <\medskipamount 
   \removelastskip \penalty 55\medskip \fi}
\def\thmm[#1]{#1}
\def\teo[#1]{#1}
\def\sttilde#1{%
\dimen2=\fontdimen5\textfont0
\setbox0=\hbox{$\mathchar"7E$}
\setbox1=\hbox{$\scriptstyle #1$}
\dimen0=\wd0
\dimen1=\wd1
\advance\dimen1 by -\dimen0
\divide\dimen1 by 2
\vbox{\offinterlineskip%
   \moveright\dimen1 \box0 \kern - \dimen2\box1}
}
\begin{document}

\title[]
{On the spectral gap of the Kac walk and other binary collision processes}

\begin{abstract}
We give a new and elementary computation of the spectral gap of the
Kac walk on $S^N$. The result is obtained as a by--product of a more
general observation which allows to reduce the analysis of the
spectral gap of an $N$--component system to that of the same system
for $N=3$. The method applies to a number of random ``binary collision'' 
processes with complete--graph structure, 
including non--homogeneous examples such as exclusion and colored exclusion processes with site disorder.   

\bigskip


  
\end{abstract}

\author[P. Caputo]{Pietro Caputo}
\address{Dip. Matematica, Universita' di Roma Tre, L.go S. Murialdo 1,
00146 Roma, Italy} \email{caputo\@@mat.uniroma3.it}

\date{July 18, 2008}

\maketitle

\thispagestyle{empty}

\section{Introduction}
The following model for energy preserving binary collisions
was introduced by M.\ Kac in \cite{Kac}.
Let $\nu$ denote the uniform probability measure on the sphere 
$$
S^{N-1} = \{\eta\in\bbR^N\,:\;\sum_{i=1}^N\eta_i^2=1\}\,,
$$
and consider the $\nu$--reversible 
Markov process on $S^{N-1}$ with infinitesimal generator given by
\be\la{genkac}
\cL f (\eta)= \frac1{2N}\sum_{i,j=1}^N \frac1{2\pi}
\int_0^{2\pi} \left[f(R_\theta^{ij}\eta) - f(\eta)\right]\,d\theta\,,
\end{equation} 
where $R_\theta^{ij}$, $i\neq j$ 
is a clockwise rotation of angle $\theta$ 
in the plane $(\eta_i,\eta_j)$. As a convention, we take 
$R^{ii}_\theta = {\rm Id}$.
In words, the associated Markov process goes as follows:
we have independent Poisson clocks of rate $1/2$ 
at each coordinate; when coordinate $i$ rings we choose uniformly at random (with replacement) another coordinate $j$; if $j\neq i$ then we perform a rotation of angle $\theta$ in the plane $(\eta_i,\eta_j)$, with $\theta$ uniform over $[0,2\pi)$, while if $i=j$ we do nothing.

Note that $-\cL$ is a non--negative, bounded self--adjoint operator on 
$L^2(\nu)$. Any constant is an eigenfunction with eigenvalue $0$ 
and the spectral gap $\l=\l(N)$ is defined as 
\be\la{gapkac}
\l(N)=\inftwo{f\in L^2(\nu):}{\nu(f)=0}\,\frac{\nu(f(-\cL) f)}{\nu(f^2)}\,,
\end{equation}
where $\nu(f)$ stands for the expectation $\int f d\nu$.
M.\ Kac conjectured that $\l(N)$ stays bounded away from $0$ as $N\to\infty$. 
This conjectured was first proved by Janvresse \cite{J}, where a powerful recursive approach due to H.T.\ Yau was used.
After that, in the beautiful paper \cite{CCL}, Carlen, Carvalho and Loss
introduced a new recursive approach which allows to actually compute
the value of $\l(N)$ for every $N$:
\be\la{ccl1}
\l(N)=\frac{N+2}{4N}\,, \quad N\geq 2\,.
\end{equation}
Around the same time, Maslen \cite{Maslen} derived formulae
for all eigenvalues of $\cL$ by means of 
harmonic analysis techniques. We refer to \cite{CCL} for further background, motivation and references on Kac's conjecture.  
Our result below shows that the proof of 
(\ref{ccl1}) can be somewhat simplified. 
In particular, we do not need any recursive analysis:
in one step we go from $\l(N)$ to $\l(3)$ and the conclusion follows by direct computations in the case $N=3$.
\begin{Th}
\la{th1}
For any $N\geq 3$:
\be
\la{bound3}
\l(N)=(3\,\l(3) - 1)\left(1-\frac2N\right) + \frac1{N}\,.
\end{equation}
In particular,  (\ref{ccl1}) follows from (\ref{bound3})
with $\l(3)=5/12$.
\end{Th} 
The proof uses a well known equivalent characterization of the spectral gap as the largest constant $\l$ such that 
the inequality 
\be
\nu\left((\cL f)^2\right) \geq \l\,\nu\left(f(-\cL )f\right)\,,
\la{equi}
\end{equation} 
holds for all $f\in L^2(\nu)$. The equivalence of (\ref{gapkac}) and
(\ref{equi}) follows from elementary spectral theory.  
A similar approach has been exploited recently in \cite{BCDP} to obtain spectral gap bounds for a class of interacting particle systems.
A proof of Theorem \ref{th1} is given at the end of the
introduction. In the following sections we shall show that 
variants of the same method can be used to obtain spectral gap estimates for several models sharing some of the features of the Kac walk. The argument turns out to be especially powerful in  
non--homogeneous cases where other known methods are harder to apply
because of the lack of permutation symmetry. In particular, for
exclusion processes with site disorder we obtain a remarkable
simplification of a spectral gap estimate proved by the author in
\cite{C1}. The latter estimate
is at the heart of recent 
results on the hydrodynamic limit of disordered lattice gases
\cite{FM,Q}. In the last section of this work we prove a new result which 
extends the spectral gap estimate to the
case of colored particles.
 
On the other hand, we point out that for some
of the physically more relevant generalizations of the Kac walk
treated in \cite{CCL}  our argument will not necessarily yield sharp
results. This is the case, for example, of the 
Kac model with momentum and energy conserving collisions discussed in
the next section.
The problem of the determination of the spectral gap for the latter model has been recently solved in \cite{CJL}, where the authors 
develop an interesting extension of the recursive scheme introduced in \cite{CCL}. 

\subsection{Proof of Theroem \ref{th1}}
We start with some preliminaries. We write 
\be\la{eij}
E_{b} f (\eta) = \frac1{2\pi}\int_0^{2\pi} f(R_\theta^{ij}\eta)\,d\theta\,,\quad b=\{i,j\}\,,
\end{equation}
for the binary average operator appearing in the definition of $\cL$.
Note that $E_{b}$ is a projection which coincides with the
$\nu$--conditional expectation given the $\si$--algebra $\cF_{b}$ generated by variables $\{\eta_\ell,\; \ell\notin b\}$.
Thus, we rewrite the Markov generator as follows:
\be
\cL f(\eta) = \frac1N \sum_{b} D_b f (\eta)\,,
\la{gen}
\end{equation}
where
the sum runs over all $\binom{N}{2}$
unordered pairs $b$ and for each 
such pair 
\be
D_b = E_b - {\rm Id}\,,\quad\;\, E_b f(\eta) = \nu(f\tc \cF_b) \,.
\la{ub}
\end{equation}
For each $b$, $D_b$ is a bounded self-adjoint operator in $L^2(\nu)$ satisfying 
$D_b^2=-D_b$.
In particular, 
\be
\nu\left(f(-\cL)g\right) = \frac1N \sum_b \nu[D_bf D_bg]\,.
\la{efg}
\end{equation}
On the other hand we have
\be
\nu\left((\cL f)^2\right) =
\frac1{N^2} \sum_{b,b'} \nu[D_b f D_{b'} f]\,.
\la{ufg1}
\end{equation}
We are going to use the expressions (\ref{efg}) and (\ref{ufg1}) in (\ref{equi}) to compute $\l(N)$.  We start with the lower bound. 
We write $b\sim b'$ when two unordered pairs have at least one common
vertex (including the case $b=b'$). Otherwise we write
$b\not\sim b'$. We observe that
if $b\not\sim b'$, then $E_b$ and $E_{b'}$ commute. Therefore,
using $D_b^2=-D_b$ and self--adjointness
\be\la{commu}
\nu[D_b f D_{b'} f] = - \nu[(D_{b'}D_b f) (D_{b'} f)] 
= \nu[(D_{b'}D_b f)^2] \geq 0\,,\quad \;b\not\sim b'\,.
\end{equation}
It follows that 
\be\la{argo1}
\nu\left((\cL f)^2\right) \geq \frac1{N^2} \sum_{b,b':\; b\sim b'} 
\nu[D_b f D_{b'} f]\,.
\end{equation}
Unordered triples $\{i,j,k\}$ of distinct vertices 
are 
denoted by $T$ (triangles). We say that $b\in T$ if $b=\{i,j\}$ and 
$i,j\in T$. Clearly, if $b\sim b'$ 
and $b\neq b'$ there is only one triangle $T$ such
 that $b,b'\in T$. 
We may therefore write
\begin{align}
\sum_{ b,b'\,:\; b\sim b'} & \nu[D_b f D_{b'}f]   = \sum_{b,b'\,:\; b\sim b'\,,\;b\neq b'} \nu[D_b f D_{b'}f] 
+ \sum_{b} \nu[(D_b f)^2]  \nonumber \\
& = \sum_{T} \sum_{b,b'\in T} \nu[D_b f D_{b'}f]  
- \sum_{T} \sum_{b\in T} \nu[(D_b f)^2]  + \sum_{b} \nu[(D_b f)^2] \,.
\la{tria}
\end{align}
Since for every $b$ there are exactly $N-2$ triangles $T$ such that
$b\in T$ we see that
\begin{align}
\sum_{ b,b'\,:\; b\sim b'} &  \nu[D_b f D_{b'}f] \nonumber\\
&\quad= \sum_{T} \sum_{b,b'\in T} \nu[D_b f D_{b'}f]  - (N-3)\sum_{b} \nu[(D_b f)^2] \,.
\la{co9}
\end{align}
Let us now apply the inequality (\ref{equi}) to a fixed triangle $T$.
Let $\cF_T$ denote the $\si$--algebra generated by
$\{\eta_\ell,\;\ell\notin T\}$. The conditional probability
$ \nu[\cdot\tc \cF_T]$ coincides with the uniform probability measure
on the sphere $S^2(t)$ in $\bbR^3$ with radius 
$$
t = \sqrt{1-\sum_{\ell\notin T}\eta_\ell^2}\,.
$$
Moreover, as noted in \cite{CCL}, it is not hard to show that 
the spectral gap of the Kac model 
does not depend on the radius of the sphere on which the walk is performed. 
Using (\ref{equi}) on each triangle $T$, we therefore have
\be
\frac13
\sum_{b,b'\in T} \nu[D_b f D_{b'} f\tc\cF_T] \geq \l(3)\sum_{b\in T} 
\nu[(D_b f)^2\tc \cF_T]\,,
\la{trgap}
\end{equation}
uniformly in $\eta\in S^{N-1}$. 
Taking $\nu$--expectation we can remove the conditioning on $\cF_T$ in
(\ref{trgap}). 
Using this in (\ref{co9}) gives
\begin{align*}
\sum_{ b,b'\,:\; b\sim b'} & \nu[D_b f D_{b'}f]\\ 
& \geq  3\,\l(3)\,(N-2) \sum_{b} \nu[(D_b f)^2]  - (N-3)\sum_{b} \nu[(D_b f)^2]  \\
& = \left((3\l(3)-1)(N-2) + 1\right)
\sum_{b} \nu[(D_b f)^2] \,.
\end{align*}
From (\ref{equi}) we conclude that $\l(N)$ is larger or equal than the right hand side of (\ref{bound3}). 

It remains to show that this bound is attained for a given $f$. 
To this end, take 
\be
f_N(\eta)=\sum_{i=1}^N \eta_i^4  + {\rm const.}
\,
\la{eigen1}
\end{equation} 
Let us first check that 
$\nu[D_b f_N D_{b'} f_N] = 0$ whenever
$b\not\sim b'$, so that (\ref{argo1}) is an equality for $f=f_N$.
To this end, note that $\nu[D_b f_N D_{b'} f_N] = \nu[f_N D_b D_{b'}
f_N] = 0$. Indeed, if $b=\{i,j\}$ and 
$b'\not\sim b$, then $D_{b'}f_N(\eta)$ depends on $\eta_i,\eta_j$ only
through $\eta_i^2+\eta_j^2 = 1- \sum_{k\notin b} \eta_k^2$ and
therefore $D_bD_{b'}f_N=0$. 

Next, we need that (\ref{trgap}) 
is an equality as well. This requires checking that for $N=3$, 
for any value of the conservation law $\sum_{i=1}^3 \eta_i^2=t>0$, 
$f_3$ is, up to additive constants (that may depend on $t$), an
eigenfunction of $-\cL$ with eigenvalue $\l(3)$. 
For the solution of this $3$--dimensional problem, as well as for the
calculation of $\l(3)=5/12$, we refer to \cite[Section 3]{CCL}.
Once the estimates (\ref{argo1}) and (\ref{trgap}) 
can be turned into identities we see that all our bounds are saturated for the function (\ref{eigen1}). This completes the proof.
%
\qed

\section{Homogeneous models}\la{hom}
The general setting can be described as follows. We consider a product space $\O=X^N$, where $X$, the single component space 
is a measurable space equipped with a probability measure 
$\mu$. 
On $\O$ we consider the product measure $\mu^N$. Elements of
$\O$ will be denoted by $\eta=(\eta_1,\dots,\eta_N)$. Next, we take a measurable function $\xi:X\to\bbR^d$, for a given $d\geq 1$, and 
we define the probability measure $\nu=\nu_{N,\o}$ on $\O$ as 
$\mu^N$ conditioned on the event 
\be \O_{N,\o}:=\left\{\eta\in\O:\;\sum_{i=1}^N\xi(\eta_i)=\o\right\}\,,
\la{cons}
\end{equation}
where  $\o\in\bbR^d$ is a given parameter.
We interpret the constraint on $ \O_{N,\o}$
as a {\em conservation law}. 

In all the examples considered below there are no difficulties in defining the conditional probability $\nu$, 
therefore we do not attempt here at a justification of this setting 
in full generality but rather refer 
to the examples for full rigor. The crucial property of $\nu$ that will be repeatedly used below is that, for any set of indices $A$, conditioned on the $\si$--algebra $\cF_A$ generated by 
variables $\eta_i$, $i\notin A$, $\nu$ becomes the $\mu$--product law over $\eta_j$, $j\in A$, conditioned on the event 
$$
\sum_{j\in A}\xi(\eta_j)=\o - \sum_{i\notin A}\xi(\eta_i)\,.
$$ 
We shall call this the {\em non--interference} property.
In analogy with (\ref{gen}) we consider the binary collision 
Markov process described by the infinitesimal generator
\be\la{genhom}
\cL f = \frac1N \sum_{b} \left[\nu(f\tc\cF_b) - f\right] \,,
\end{equation}
where
the sum runs over all $\binom{N}{2}$
unordered pairs $b=\{i,j\}$ 
and 
$\nu[ f\tc\cF_b] $ is the $\nu$--conditional expectation of $f$ given
the variables $\eta_\ell,\; \ell\notin b$. 
This defines a bounded self--adjoint operator on $L^2(\nu)$.
Setting, as before, $D_b =\nu[\cdot\tc \cF_b] - {\rm Id}$, we see that 
$D_b^2=-D_b$ and the operator $\cL$ satisfies
(\ref{efg}).

In principle, our arguments could be extended to more general
processes. For instance, one can consider binary ``collisions'' which are not
given by simple averages as in (\ref{genhom}) but by other mechanisms
which still preserve reversibility (see e.g.\ the non uniform models
considered in \cite{CCL}), or one could look at 
``collisions'' which can involve 
more than two components at a time. However, we
shall not investigate such extensions.

Returning to our model (\ref{genhom}), 
the spectral gap is defined as in (\ref{gapkac}). 
Note that, by definition,  
$\l(2)=\frac12$ always, since for $N=2$ there is only one pair $b$ and $\nu((D_b f)^2) = \nu(f^2)$ 
for any $f$ such that $\nu(f)=0$. 
Here we shall assume that 
$\l(3)$ {\em is independent of the choice of} $\o$. 
This is a strong assumption which holds only for special choices of 
the model. However, it does hold in the examples considered below.  
More general models are treated in the next section.
\begin{Th}
\la{th2}
Suppose $\l(3)$ is independent of $\o$.
Then, for any $N\geq 2$:
\be
\la{bar3}
\l(N)\geq (3\,\l(3) - 1)\left(1-\frac2N\right) + \frac1{N}\,.
\end{equation}
If, in addition, 
there exists $\varphi:X\to\bbR$ such that the function  
\be\la{sumform}
f_3(\eta_1,\eta_2,\eta_3)=\sum_{i=1}^3 \varphi(\eta_i)\,, 
\end{equation}
satisfies, for $N=3$, $\cL f_3 = -\l(3) f_3 + {\rm const.}$, 
regardless of the value of 
$\sum_{i=1}^3\xi(\eta_i)$ (although the constant may depend on this value), 
then (\ref{bar3}) can be turned into an identity for each $N\geq 2$. 
\end{Th} 
\proof 
We repeat the steps of the proof of Theorem \ref{th1}. We start from (\ref{ufg1}) and arrive at (\ref{argo1}) with the same commutation property used in (\ref{commu}). Indeed, this is a simple consequence of the non--interference property. 
The latter property also implies that the conditional probability $\nu(\cdot\tc\cF_T)$ is nothing but $\nu_{3,\o'}$ with $\o'=\o-\sum_{\ell\notin T} \xi(\eta_\ell)$.
Since $\l(3)$ is independent of the value $\o$ in the conservation 
law we may repeat the argument leading to (\ref{trgap}).
This proves the lower bound (\ref{bar3}).
As for the reverse direction, again the arguments given in the proof
of Theorem \ref{th1} can be repeated line by line. \qed

\bigskip
Next, we examine some examples to which the theorem applies.

\subsection{Kac model} The model discussed in the introduction can be
seen as a special case of our general setting, so that Theorem
\ref{th1} becomes a special case of Theorem \ref{th2}. 
Here $X=\bbR$, $\mu$ is the centered Gaussian measure with variance $\si^2>0$ and we take 
$\xi(\eta_i)=\eta_i^2$ (with $d=1$). Then, for every $\o>0$,  
$\nu_{N,\o}$ is the uniform probability measure on the sphere of radius
$\sqrt \o$. Clearly, the choice of $\si^2>0$ is uninfluential in the determination of $\nu_{N,\o}$. 
As we have seen in the introduction, 
this model satisfies the two main assumptions in Theorem \ref{th2}. 

\subsection{``Flat'' Kac model}  
Here $X=\bbR^+$, $\mu$ is the exponential law with parameter $\g>0$. We take $\xi(\eta_i)=\eta_i$ (with $d=1$). Then, for every $\o>0$,
independently of the choice of $\g>0$,  
$\nu_{N,\o}$ is the uniform probability measure on the simplex 
$\O_{N,\o}$. The binary collision process (\ref{genhom}) associated to this
setting does not appear explicitly in the literature, so we shall give
more details on the computation of $\l(N)$ in this case. 
Let us first 
check that $\l(N)$ is independent of $\o$, for any $N$. 
This is a consequence of the fact that
$\cL$ commutes with the unitary change of scale from $\O_{N,\o}$ to
$\O_{N,\o'}$, for any $\o,\o'>0$. Indeed, $\nu_{N,\o'}$ is the image
of $\nu_{N,\o}$ under the map $\cT:\eta \to \o'\eta/\o$ and if
$f_\cT(\eta)=f(\cT\eta)$, then 
\be\la{comma}
\nu_{N,\o}(f_\cT\tc\cF_b)(\eta) = \nu_{N,\o'}(f\tc\cF_b)(\cT\eta)\,,
\end{equation} 
for all $\eta\in\O_{N,\o}$ and for all pairs $b$. 

Next, we shall prove that $\l(3)=\frac49$ and that the eigenfunction
for $N=3$ is given by 
\be
\la{sumflat}
f_3(\eta_1,\eta_2,\eta_3) =  \sum_{i=1}^3 \eta_i^2 + {\rm const.}
\end{equation}
for any value of the constraint $\sum_{i=1}^3\eta_i=\o$
(of course, the constant will be given by $-3\nu_{3,\o}(\eta_1^2)$, since we must have $\nu_{3,\o}(f)=0$).
From these facts and Theorem \ref{th2} we therefore obtain, for all $N\geq 2$: 
\be\la{gapflat}
\l(N) = \frac{N+1}{3N}\,.
\end{equation}
To solve the $3$--dimensional problem, we observe that when $N=3$, then $\cL+1$ coincides with the {\em average operator} $P$ introduced in \cite{CCL}. Therefore we can apply the general analysis of
Section 2 in \cite{CCL} or equivalently that of Theorem 4.1 in \cite{C2}. The outcome is that 
\be
\l(3)\geq \frac13\,\min\{2+\mu_1\,,\,2-2\mu_2\}\,,
\la{3k}
\end{equation}
where the parameters $\mu_1,\mu_2$ are given by 
\be\la{mus}
\mu_1=\inf_{\varphi}
\nu(\varphi(\eta_1)\varphi(\eta_2))
\,,\quad\; 
\mu_2=\sup_{\varphi}
\nu(\varphi(\eta_1)\varphi(\eta_2))
\end{equation}
with $\varphi$ chosen among all functions $\varphi: X\to \bbR$ satisfying 
$\nu( \varphi(\eta_1)^2)=1$ and $\nu(\varphi(\eta_1))=0$.
Here $\nu$ stands for $\nu_{3,\o}$, but we have removed the subscripts
for simplicity. As in (\ref{comma}) one checks that 
the parameters $\mu_1,\mu_2$ do not depend on $\o$. 
Write $\cK\varphi(\a) = \nu[\varphi(\eta_2)\tc\eta_1=\a]$,
$\a\geq 0$. This defines a self--adjoint Markov operator on
$L^2(\nu_1)$, where $\nu_1$ is the marginal on $\eta_1$ of
$\nu$. 
In particular,
the spectrum ${\rm Sp}(\cK)$ of $\cK$ 
contains $1$ (with eigenspace given by the
constants). Then $\mu_1,\mu_2$ are, respectively, the smallest and the
largest value in ${\rm Sp}(\cK)\setminus\{1\}$, as we see by writing 
$\nu(\varphi(\eta_1)\varphi(\eta_2)) =
\nu[\varphi(\eta_1)\cK\varphi(\eta_1)]$. This is now a
one--dimensional problem and $\mu_1,\mu_2$ can be computed as
follows. To fix ideas we use the value $\o=1$ for the conservation law 
$\eta_1+\eta_2+\eta_3$. In this case $\nu_1$ is the law on 
$[0,1]$ with density $2(1-\eta_1)$.
 Moreover,
$$
\cK\varphi(\eta_1) = \frac1{1-\eta_1}\int_0^{1-\eta_1}\varphi(\eta_2)d\eta_2\,, \quad \eta_1\in[0,1)\,.
$$
($\cK\varphi(1)=\varphi(0)$). 
In particular, $\varphi_1(\a)=\a - \frac13$ is an eigenfunction of
$\cK$ with eigenvalue $-1/2$. 
Moreover, $\cK$ preserves the degree of polynomials so that if $Q_n$
denotes the space of all polynomials of degree $d\leq n$ we have $\cK Q_n\subset
Q_n$. By induction we see that for each $n\geq 1$ the
polynomial
$\a^n + q_{n-1}(\a)$, for a suitable 
$q_{n-1}\in Q_{n-1}$, is an eigenfunction with eigenvalue  
$\mu_n=\frac{(-1)^{n}}{n+1}$, and it is orthogonal to
$Q_{n-1}$ in $L^2(\nu_1)$. Since the union of $Q_n$, $n\geq 1$, is dense 
in $L^2(\nu_1)$ 
this shows that there is a complete orthonormal set of eigenfunctions
$\varphi_n$, where $\varphi_n$ is a polynomial of degree $n$ with
eigenvalue $\mu_n$ and ${\rm Sp}(\cK)=\{\mu_n\,,\;n=0,1,\dots\}$.
Therefore we can take 
$\mu_1=-\frac12$ and 
$\mu_2 = \frac13$ in the formula (\ref{3k}). We conclude 
that $\l(3)\geq \frac49$. 

To end the proof we take $f=\eta_1^2+\eta_2^2+\eta_3^2$
and, using $\nu[\eta_1^2\tc\eta_2] = \frac13(\eta_2^2 - 2\eta_2 + 1)$ we compute
$$
\cL f (\eta) = -\frac49 \,f(\eta) + {\rm const.}\,
$$
Thus, $\l(3)=\frac49$ and the eigenfunction is given by (\ref{sumflat}).
Clearly, the unitary change of scale $\cT$ introduced above does not alter the form of the eigenfunction so that all the hypothesis of Theorem \ref{th2} apply and (\ref{gapflat}) follows.

\subsection{Momentum and energy conserving collision model}. Here
$X=\bbR^3$, and $\mu$ is a centered $3$--dimensional Gaussian law
$N(0,C)$ with covariance matrix $C$ given by a multiple of the
identity. Each coordinate $\eta_i$ is a 
$3$--dimensional velocity vector $\eta^\a_i$, $\a=1,2,3$. We have $d=4$ conservation laws, with 
$\xi^\a(\eta_i)=\eta^\a_i$, $\a=1,2,3$ (momentum conservation)
and $\xi^4(\eta_i) = |\eta_i|^2=\sum_{\a=1}^3(\eta^\a_i)^2$
(energy conservation).
For any $\o=(\o^1,\dots,\o^4)\in\bbR^4$ with $\o^4 > 0$, $\nu_{N,\o}$
is the uniform probability measure on the manifold $\O_{N,\o}$
(whenever $\O_{N,\o}\neq\emptyset$). 
We refer to \cite{CCL} for an explicit description of the 
probability measure and its main properties. It is still the case that $\l(3)$ is independent of the conservation law, see \cite{CCL}. In particular, the lower bound (\ref{bar3}) holds in this case. However, a computation of $\l(3)$ for this model shows that $\l(3)=\frac13$ (see (\ref{gapbol}) below) and therefore the estimate becomes $\l(N)\geq \frac1N$ which is rather poor. 
Indeed, it was shown in \cite{CCL} that the 
the spectral gap is bounded away from zero uniformly in $N$: 
\be\la{inf}
\inf_{N\geq 2} \l(N)>0\,.
\end{equation}
Recently, by a very deep analysis of the Jacobi polynomials naturally associated to this model Carlen, Jeronimo and Loss \cite{CJL} 
succeeded in computing $\l(N)$ exactly for every $N$:
\be
\la{gapbol}
\l(N)=\frac13\,,\quad\; N\geq 3\,.
\end{equation} 
This shows our approach is too rough here. 
As we know from Theorem \ref{th2} the loss must come from the lack of the second property required in that theorem. It was shown in \cite{CJL} that the eigenspace of $\l(N)$, for the choice $\o=(0,0,0,1)$, is spanned by
the functions  
$$
f_{N,\a}(\eta)=\sum_{i=1}^N |\eta_i|^2\eta^\a_i\,,\quad\a=1,2,3\,.
$$
One cannot expect that the 
change of scale from $\o$ to $\o'$ 
transforms a linear combination of $f_{N,\a}$'s into itself 
(up to multiplicative and additive constants), and 
the second property in Theorem \ref{th2} must fail here.

Let us show that our approach can nevertheless be used to prove the
weaker result (\ref{inf}) without any recursive analysis. Namely, we
prove that if $\l(4)>1/4$ then (\ref{inf}) holds.
The choice of triangles in the proof of Theorem \ref{th1} and Theorem \ref{th2} can be replaced by the choice of larger cliques
(i.e.\ complete subgraphs) of the original complete graph. 
Namely, if in (\ref{tria}) we sum over cliques with $4$ vertices instead of triangles we shall obtain the bound, for $N\geq 4$:
\be\la{l4}
\l(N)\geq (4\l(4) - 1)\left(\frac12-\frac1N \right) + \frac1N\,,
\end{equation}
instead of (\ref{bar3}). To see this, set for simplicity 
$a_{b,b'}:=\nu[D_bf D_{b'}f]$, for a given $f\in L^2(\nu)$, 
and recall that $N^2\nu[(\cL f)^2] = \sum_{b,b'}a_{b,b'}$. Denote by $Q$ the cliques of $4$ vertices and note that: for every $b$ there are $\frac12(N-2)(N-3)$  $Q$'s such that $b\in Q$; for any $b\neq b'$ with $b\sim b'$ there are $(N-3)$ $Q$'s 
such that $b,b'\in Q$; for any $b, b'$ with $b\not\sim b'$ there is
only one $Q$ 
such that $b,b'\in Q$. Then 
\begin{align*}
\sum_{b,b'\,:\, b\neq b',\, b\sim b'}  a_{b,b'} & = \frac1{N-3}\sum_{Q}
\sum_{b,b'\in Q\,:\, b\neq b',\, b\sim b'}  a_{b,b'} \\ & = 
\frac1{N-3}\sum_{Q}
\left\{\sum_{ b, b'\in Q}  a_{b,b'} - \sum_{ b, b'\in Q\,:\; b\not\sim b'}  a_{b,b'} - \sum_{ b\in Q}  a_{b,b}\right\}\\ &
 \geq \frac1{N-3}\sum_{Q}(4\l(4) - 1) \sum_{ b\in Q}  a_{b,b}
- \frac1{N-3}\sum_{Q}\sum_{ b, b'\in Q\,:\; b\not\sim b'}  a_{b,b'} \\&
= \frac12(N-2)(4\l(4) - 1) \sum_{ b}  a_{b,b} - \frac1{N-3}\sum_{ b, b'\,:\; b\not\sim b'}  a_{b,b'} \,.
 \end{align*}
 Therefore 
  \begin{align*}
\sum_{ b,b'}  a_{b,b'} 
 &= \sum_{ b\neq b'\,:\; b\sim b'}  a_{b,b'}  
 + \sum_{ b, b'\,:\; b\not\sim b'}  a_{b,b'} +  \sum_{ b}  a_{b,b}\\
& \geq \frac12 [(N-2)(4\l(4)-1) + 2]  \sum_{ b}  a_{b,b} 
+ \left(1-\frac1{N-3}\right)
 \sum_{ b, b'\,:\; b\not\sim b'}  a_{b,b'}
 \\ & \geq  \frac12[(N-2)(4\l(4)-1) + 2] \sum_{ b}  a_{b,b}\,,
 \end{align*}
 where, in the last line, we have used (\ref{commu}).
 Since $-N\nu(f\cL f) = \sum_{ b}  a_{b,b}$, this proves the claim (\ref{l4}). Note that this argument applies in the more general setting of Theorem \ref{th2}, and similar computations can in principle be carried out for any choice of cliques of $k<N$ vertices. 
Therefore, to prove (\ref{inf}) 
it suffices to prove $\l(4)>\frac14$. With the value $\l(4)=\frac13$ from (\ref{gapbol}) this gives $\l(N)\geq \frac16 + \frac2{3N}$, for 
all $N\geq 4$. 


\section{Non--homogeneous models}\la{nonhom}
We generalize the setting introduced above as follows. As before we
take $\O=X^N$ but now each copy of $X$ will be equipped with 
possibly distinct probability measures $\mu_i$, $i=1,\dots,N$. Again
we consider conservation laws as in (\ref{cons}), associated to a
given function $\xi$ on $X$ and a given parameter $\o$, and the
probability $\nu$ given by the product $\mu_1\times\cdots\times\mu_N$ conditioned on
$\eta\in\O_{N,\o}$. Note that the non--interference property still
holds for this setting.  
 The binary collision process
is defined as in (\ref{genhom}). Again, $-\cL$ is a non--negative self--adjoint operator on $L^2(\nu)$ and we may define its spectral gap just as in (\ref{gapkac}). This time, however, to keep track of the conservation law we shall write $\l(N,\o)$ instead of just $\l(N)$. Note that $\l(2,\o)=\frac12$ always. We define
\be
\bar\l(N) = \inf_{\o} \l(N,\o)\,,
\la{infl}
\end{equation}
where the infimum ranges over the set of admissible values of the parameter 
$\o$. This set depends on the choice of the model and, as usual, we
refer to the examples for fully rigorous formulations of the results. 
As a convention we may set $\l(N,\o)=+\infty$ if $\o$ is such that the
measure $\nu$ becomes a Dirac delta.  
Thanks to the non--interference property of $\nu$ 
there is no difficulty in repeating the 
previous arguments to prove the following estimate.
\begin{Th}
\la{th3}
For any $N\geq 2$ and any $\o$:
\be
\la{barbar3}
\bar\l(N)\geq (3\,\bar\l(3) - 1)\left(1-\frac2N\right) + \frac1{N}\,.
\end{equation}
\end{Th} 

Let us investigate some specific models.

\subsection{Non--homogeneous Kac models}
Consider the following non--homogeneous version of the ``flat'' Kac model introduced in Section \ref{hom}. Take $X=\bbR_+$ and $\mu_i$ the probability on $\bbR_+$ with density 
$$
\frac1{z_i}\exp{(-\eta_i + b_i(\eta_i))}\,,$$
where $b_i$  are bounded 
measurable functions and $z_i=\int_0^\infty dx\exp{(-x + b_i(x))}$ 
is the normalizing constant.
We set $B:=\sup_i|b_i|_\infty$. For this model, any value $\o>0$ of
the conservation law is allowed in the definition (\ref{infl}) of
$\bar \l(N)$. 
We claim that 
for every $\e<\frac13$ there is $\d>0$ such that 
\be\la{infnonh}
\liminf_{N\to\infty} \bar\l(N) \geq \e\,,
\end{equation} provided $B<\d$.
Thanks to Theorem \ref{th3}, to prove (\ref{infnonh}) 
it suffices to show that  
$\l(3,\o)\geq \frac4{9}(1-\e(B))$ with some $\e(B)\to 0$
as $B\to 0$, uniformly in $\o>0$, the value of the conservation law $\sum_{i=1}^3\eta_i = \o$. 
To this end we shall use a standard perturbation argument. 
Let $\nu$ denote the measure $\mu_1\times \mu_2\times \mu_3(\cdot\tc \sum_{i=1}^3\eta_i = \o)$  and call $\nu_0$ the same measure in the case $b_1=b_2=b_3=0$. Thus, for any bounded measurable function $g$ we have
\be\la{nuf}
\nu(g)= \frac1{Z_\o}\,\int_0^\o d\eta_1\int_0^{\o-\eta_1}d\eta_2\,
g(\eta_1,\eta_2,\o-\eta_1-\eta_2) \,u(\eta_1,\eta_2)\,,
\end{equation}
where $u(x,y)=\nep{-b_1(x)-b_2(y) - b_3(\o-x-y)}$ and $Z_\o=\int_0^\o dx\int_0^{\o-x}dy\,
u(x,y)$. Using 
$$\nep{-3B}\leq u(\eta_1,\eta_2)\leq \nep{3B}\,,$$
it is easily seen that, for any bounded $f$ we have the bound between variances   
\be\la{flat1}
\var_{\nu}(f)\leq \nep{6B}\var_{\nu_0}(f)
\,.
\end{equation}
(Use $g= (f-\nu_0(f))^2$ in (\ref{nuf}) and the fact that 
$\nu(g)\geq \var_{\nu}(f)$).
The same reasoning shows that 
\be\la{flat2}
\nu(f(-\cL) f)\geq \nep{-10B}\,\nu_0(f(-\cL_0) f)\,,
\end{equation}
where $\cL_0$ 
is the generator corresponding to the choice $b_i\equiv 0$.
Indeed,
\begin{align*}
\nu(f(-\cL) f) &= \frac13 \sum_{i=1}^3\nu\left[(\nu[f\tc\eta_i] -
  f)^2\right]\\
&=
\frac13 \sum_{i=1}^3\nu\left(\var_{\nu}(f\tc\eta_i)\right)\,,
\end{align*}
(where $\var_{\nu}(f\tc\eta_i)$ denotes the variance of $f$ w.r.t.\
$\nu(\cdot\tc\eta_i)$). For each $i=1,2,3$ we have (as above) 
\be\la{flat3}
\var_{\nu}(f\tc\eta_i)\geq \nep{-4B} \var_{\nu_0}(f\tc\eta_i)\,,
\end{equation} 
uniformly in $\eta$.
A further comparison gives 
$\nu\left(\var_{\nu}(f\tc\eta_i)\right)\geq \nep{-10B}
\nu_0\left(\var_{\nu_0}(f\tc\eta_i)\right)$ which implies (\ref{flat2}).
 
Recall that the spectral gap for $\cL_0$ is equal to $4/9$ regardless of
the value of $\o>0$. The previous estimates therefore imply that
$$
\l(3,\o)\geq \frac49\,\nep{-16 B}\,.
$$ 
This proves our claim with $\e(B) = 1-\nep{-16 B}$, from which (\ref{infnonh}) follows.

\smallskip
The same argument can be used to produce uniform lower bounds on the gap of non--homogeneous versions of the Kac walk on $S^{N-1}$
and of the momentum and energy conserving collision model,  under the assumption of  small perturbations. In the latter model we need the argument in (\ref{l4}) to obtain a uniform estimate $\inf_N\bar\l(N)>0$.
 

\subsection{Non--uniform random permutations}
We take $X=\{1,\dots,N\}$ and, for each $i=1,\dots,N$ we consider a
probability $\mu_i$ on $X$ given by
$\mu_i(\eta_i=j)=\frac{\nep{-b_i(j)}}{Z_i}$, where $b_i:X\to\bbR$ are
  bounded functions, and $Z_i=\sum_{j=1}^N\nep{-b_i(j)}$. To model permutations we use $N$ conservation
  laws that will force all components of $\eta$ to have distinct
  values: set $d=N$ and define $\xi=(\xi^j)_{j=1}^N$ with
  $\xi^j(\eta_i) = 1_{\{\eta_i=j\}}$. Fixing $\o=(1,1,\dots,1)$ we see
    that the set $\O_{N,\o}$ coincides with the set of $N!$ permutations of
    $N$ letters. 
We define the probability $\nu$ as usual by
    $\mu_1\times\cdots\mu_N(\cdot\tc\eta\in\O_{N,\o})$. Note that if
    the bias functions $b_i$ are all $0$ then $\nu$ is simply the uniform
  probability measure over permutations. 

The binary collision is now a {\em random transposition} process. 
Note that only the value $\o=(1,1,\dots,1)$ is considered for the
conservation law so that $\bar\l(N) = \l(N,\o)$ for this model (there
is no real infimum in (\ref{infl}) here). 
In the uniform case ($b_i\equiv 0$) a simple variant of the
    argument of Theorem \ref{th2} proves that $\l(N)=\frac12$ for
    every $N\geq 2$. 

While relaxation to equilibrium for the uniform
    case is well known, the non--uniform case is
    certainly less understood. Here we can show that if
    $B=\sup_i|b_i|_\infty$ is sufficiently small we have 
$\inf_N\bar \l(N) > 0$. More precisely, for every $\e<\frac12$ there exists
$\d>0$ such that $B<\d$ implies $$\liminf_{N\to\infty}\bar\l(N)\geq \e\,,$$ for all $N\geq
2$. To prove this we use exactly the same argument we have used to prove
(\ref{infnonh}). In particular, it suffices to show that
$\bar\l(3)\to\frac12$ as $B\to 0$. 
The same perturbation argument will yield the desired estimate. To avoid
repetitions we leave the details to the reader. 

\subsection{Exclusion processes with site disorder}
Here we consider a non--homogeneous version of what is sometimes called the Bernoulli--Laplace process. 
The inhomogeneous 
distribution models site impurities or site disorder.  We take $X=\{0,1\}$, $\mu_i$ is a Bernoulli law with parameter $p_i\in(0,1)$, i.e.\ $\mu_i (\eta_i=1) = p_i$
and $\mu_i (\eta_i=0) = 1-p_i$. The value of $\eta_i$ is interpreted as the presence ($\eta_i=1$) or absence ($\eta_i=0$) of a particle at the vertex $i$. The function $\xi$ is given by $\xi(\eta_i)=\eta_i$ so that 
for any integer $\o\in\{0,1,\dots,N\}$, the set $\O_{N,\o}$ denotes the configurations of $\o$ particles over $N$ vertices.  
The binary collision process (\ref{genhom}) becomes nothing but  the
well known exclusion process on the complete graph $\{1,\dots,N\}$. This can be seen as follows.
Given a pair $b=\{i,j\}$ and a configuration $\eta\in\O_{N,\o}$, write
$\o_{ij}=\o-\sum_{\ell\notin b}\eta_\ell$. Clearly, $\o_{ij}\in\{0,1,2\}$. 
Observe that, if $\o_{ij}=1$ then $\nu(f\tc\cF_b)(\eta) $ is given by
$$
 \frac{p_i(1-p_j)}{p_i(1-p_j)+p_j(1-p_i)} f(\eta;1,0)
+ \frac{p_j(1-p_i)}{p_i(1-p_j)+p_j(1-p_i)} f(\eta;0,1)
$$
where, for simplicity we write explicitly the $i$-th and $j$-th entries in
$f(\eta)=f(\eta;\eta_i,\eta_j)$. On the other hand in the case
$\o_{ij}\in\{0,2\}$ we have $ \nu(f\tc\cF_b)(\eta)  = f(\eta)$.
Setting 
\be\la{ratesoo}
c_{b}(\eta) = \frac{p_i(1-p_j)\eta_j(1-\eta_i)}{p_i(1-p_j)+p_j(1-p_i)}
+ \frac{p_j(1-p_i)\eta_i(1-\eta_j)}{p_i(1-p_j)+p_j(1-p_i)}\,,
\end{equation}
we therefore obtain, for any $ \eta\in\O_{N,\o}$,
\be\la{rateso}
\nu(f\tc\cF_b)(\eta) - f(\eta)= c_b(\eta)\left(f(\eta^b) - f(\eta)\right)\,,
\end{equation}
where $\eta^b$ denotes the configuration 
in which $\eta_i$ and $\eta_j$ have been exchanged. 
From (\ref{rateso}) wee see that $\cL$ has the familiar form of the 
exclusion process.

If we proceed by perturbative arguments (as in the previous
two subsections) we would be able to prove a
result of the form: if $p_i$ are (uniformly) sufficiently close to
$\frac12$
then we have a uniform bound from below on the spectral gap. 
However, we shall prove here a much stronger result. 
We assume there exists $\e>0$ such that the parameters $p_i$ satisfy
\be
\la{pi}
\e\leq p_i\leq 1-\e\,,\quad i=1,\dots,N\,.
\end{equation}
The minimal spectral gap $\bar\l(N)$ is defined as usual by
(\ref{infl}), where the infimum ranges over all
$\o\in\{0,1,\dots,N\}$, with the convention that 
$\l(N,0)=\l(N,N)=\infty$.
Under the same assumption the following theorem was proved in
\cite{C1} by means of rather technical local limit theorem
estimates. It is surprising that the simple argument of Theorem
\ref{th3} allows a straightforward proof. The uniform spectral gap
bound below is an important step in the recent works \cite{FM,Q} establishing hydrodynamic limits for exclusion processes with disorder. 

\begin{Th}\la{BEgg}
Assume (\ref{pi}) for some $\e>0$. Then, there exists $c_\e>0$ 
such that for all $N\geq 2$ 
\be\la{exc1}
\bar\l(N)\geq c_\e\,.
\end{equation}
\end{Th}
\proof
Thanks to Theorem \ref{th3}, all we have to do is prove that 
\begin{equation}\la{exc2}
\bar\l(3)\geq \frac13+ c_\e\,,
\end{equation}
 for some $c_\e>0$.
We fix three vertices $i=1,2,3$, with their occupation probabilities $p_i$
satisfying (\ref{pi}) and with $\sum_{i=1}^3\eta_i = \o$.
We may assume that $\o=1$, i.e.\ 
there is one particle. Indeed if
there are two we may look at occupied vertices as empty and
vice-versa, if there is none (or three) the measure is a Dirac delta
and 
by our convention (see discussion after (\ref{infl})) 
the estimate $\l(3,\o)\geq  \frac13+ c_\e$ becomes obvious. 
Since there is one particle we shall call $x$, respectively $y$,
the probability that the particle is at $i=1$, respectively $i=2$.
We set $z=1-x-y$ for the probability that the particle is at $i=3$. 
Note that, thanks to (\ref{pi}), $x,y,z$ are all bounded away from $0$ and $1$. For instance,
$$
x= \frac{p_1(1-p_2)(1-p_3)}{p_1(1-p_2)(1-p_3)+(1-p_1)p_2(1-p_3)+(1-p_1)(1-p_2)p_3}\,.
$$
It is easily checked that the process generated by $\cL$ on our three sites becomes 
a $3$--state Markov chain with the $3\times 3$ transition
matrix $P=\cL+{\rm Id}$ given by
\be
P = \frac13 \,\left( \begin{array}{ccc}
1+
\frac{x}{x+y}
+
\frac{x}{x+z} 
& \frac{y}{x+y} & \frac{z}{x+z}  \\
\frac{x}{x+y} & 1+
\frac{y}{x+y}
+
\frac{y}{y+z}  & 
\frac{z}{y+z}  \\
\frac{x}{x+z}  & \frac{y}{y+z}  & 1+
\frac{z}{x+z}
+
\frac{z}{y+z} 
 \end{array} \right)\,.
\la{P}
\end{equation}
We need to estimate the eigenvalues of $P$. Clearly, one
eigenvalue is $1$. From (\ref{P}) we see that ${\rm Tr}(P)=2$.
Therefore the other two
eigenvalues must satisfy $\l_1={\rm Tr}(P)-1-\l_2=1-\l_2$.
Note that $$\l(3,\o)=\min\{1-\l_1,1-\l_2\}\,.$$ 
To estimate $\l_i, i=1,2$ we compute the determinant of $P$. 
The next lemma shows that $\det(P)>\frac29$. Therefore, for both
$i=1,2$ we have 
$\l_i(1-\l_i) > \frac29$. This implies $|\l_i - \frac12| < \frac16$. 
In particular, it shows that $\l_i < \frac23$. In conclusion:
$\bar\l(3)>\frac13$ as claimed. (Note that $\l_i=\frac12$ would be the value in the homogeneous case $p_i\equiv\frac12$.)
\begin{Le}
\la{detp}
\be
\det(P) = \frac29\left( 1 + \frac{xyz}{(1-x)(1-y)(1-z)}\right)
\la{detp1}
\end{equation}
\end{Le}
\proof 
Set $\G=3 P$. Also, set $$
\d = \frac{x}{x+z}\,,\quad 
\b = \frac{x}{x+y}\,,\quad
\g = \frac{y}{y+z}\,.
$$
From (\ref{P}) we compute
\begin{align*}
\det(\G)& = \left(1+ \b + \d\right)\left(6 -
  3\,\b - 2\,\d  +
  \b\,\,\d
- \g\,\,\d + \g\,\,\b
\right) \\
&+ \left(\b -1 \right)
\left(3\,\b -  \b\,\,\d 
- \g\,\,\b - \d + \g\,\,\d
\right)\\
&+ \left(1-\d \right)
\left(\b\,\,\d  
- \g\,\,\d + 
\g\,\,\b - 2\,\d \right)\,.
\end{align*}
Simplifying we arrive at 
$$
\det(\G)= 6 + 3\,\d - 3\,\b\,\d -   3\,\g\,\d + 3\,\b\,\g
\,.
$$
Rewriting this in terms of the probabilities $x,y,z$ we 
obtain 
\be
\det(\G) = 6 + \frac{6\,x\,y\,z}{(1-x)\,(1-y)\,(1-z)}\,.
\la{detp3}
\end{equation}
This implies (\ref{detp1}). \qed

\subsection{Colored exclusion processes with site disorder}
A natural generalization of the previous model is a system
where particles can be of several different kinds - or colors.
Namely, suppose there are $m$ colors and let each 
particle be painted with one of the available colors. 
Configurations of colored particles are denoted by
$\eta\in\O:=\{0,1,\dots,m\}^N$ with the
interpretation that $\eta_i=0$ means that $i$ is empty, while
$\eta_i=k$, $k\in\{1,\dots,m\}$ means that $i$ is occupied by a
particle with color $k$. Thus, the single state space $X$ is 
$\{0,\dots,m \}$. The conservation laws are given by 
$\xi^k(\eta_i)=1_{\{\eta_i = k\}}$, 
$k=1,\dots,m$ and the vector $\o=(\o_1,\dots,\o_m)$, where 
$\o_i$ are non--negative integers such that $\sum_{k=1}^m\o_k \leq N$. Thus, the set $\O_{N,\o}$ denotes the configurations of particles over $N$ vertices, 
with $\o_k$ particles with color $k$. 
We shall use the notation 
$\psi_i=1_{\{\eta_i\geq 1\}}$ so that the variables
$\psi\in\{0,1\}^N$ denote the configuration of occupied sites.
Let $p_i\,,\;i=1,\dots N$ be given parameters satisfying (\ref{pi}). 
Let $\mu_i$ denote the probability on $X$ such that
$$
\mu_i(\eta_i=k)=\frac1{Z_i}\left(p_i1_{\{k \geq 1\}} + (1-p_i)1_{\{k = 0\}}\right)\,,
$$
where $Z_i=(m-1)p_i +1$ is the normalization.
In particular, w.r.t.\ $\mu_i$ the occupation variable $\psi_i$ is a Bernoulli random variable with parameter $mp_i/Z_i$. We call, as usual, $\nu$ the product
$\mu_1\times\cdots\mu_N$ conditioned on $\O_{N,\o}$. To avoid
degenerate cases we take $1\leq \sum_{k=1}^m\o_k \leq N-1$.
We are interested in the following exclusion--type dynamics.
For any configuration $\eta$ and any edge $b=\{i,j\}$ we write 
$\eta^b$  
for the configuration where the variables $\eta_i,\eta_j$ have been exchanged. 
For every $\g\geq 0$ we define the Markov generator by
\be
\cL_\g f(\eta) = \frac1{N}\sum_{b}c^\g_b(\eta)
\left(f(\eta^b)-f(\eta)\right)\,,
\la{genera}
\end{equation}
where the rates are given, for $b=(i,j)$, by
\be
c^\g_{b}(\eta) = c_b(\psi) + \frac\g{2}\,1_{\{\psi_i=\psi_j\}}
\la{rates}
\end{equation}
where $c_b$ are the functions defined in (\ref{ratesoo}), but now evaluated at the occupation variables $\psi=\psi(\eta)$ defined above.
It is easily checked that the rates (\ref{rates}) are reversible w.r.t.\
$\nu$: for any $\eta$ and any $b$
\be
\nu(\eta)c^\g_b(\eta) = \nu(\eta^b)c^\g_b(\eta^b)\,.
\la{detba} 
\end{equation} 
The latter statement is equivalent to self--adjointness of $\cL_\g$
in $L^2(\nu)$. Moreover, 
\be\la{diri}
-\nu(f\cL_\g f) = \frac12\sum_{b}\nu\left[c^\g_b\,(f^b-f)^2\right]\,,
\end{equation}
where $f^b(\eta):=f(\eta^b)$. 
If $1\leq \sum_{k=1}^m\o_k \leq N-1$ the Markov chain 
generated by $\cL_\g$ is irreducible. We shall consider the cases 
$\g=0$ and $\g=1$. 
The difference between $\cL_0$ and $\cL_1$ is that in $\cL_1$ we have
added the possibility of ``stirring'' between particles, i.e.\
exchange of positions of particles of different colors. The case
$\g=1$ coincides then with the usual binary collision dynamics given
by local averages (\ref{genhom}). 
On the other hand, the case 
$\g=0$ is a true {\em exclusion} process, with particles jumping only
to empty sites. The addition of stirring can result in a faster
relaxation to the equilibrium distribution $\nu$, or equivalently in a
larger spectral gap. Note, however, that for the case $m=1$ there is
no difference: in this case $\cL_\g=\cL_0$ for all $\g$. The case
$m=1$, of course, is the one analyzed in  
Theorem \ref{BEgg}. From now on we take $m>1$.

We denote by $\l_\g(N,\o)$ the spectral gap
of the generator $\cL_\g$. The following theorem shows that 
in the case $\g=1$ we have a uniform lower bound $\l_1(N,\o)\geq c_\e$ as in Theorem \ref{BEgg} and, in the case $\g=0$ we
have $ \l_0(N,\o) \geq c_\e(1-\r)$, where $\r$ is the global density:
\be
\r=\frac1N\sum_{k=1}^m \o_k\,.
\la{rho}
\end{equation}
The slow--down in the limit $\r\to 1$ is natural in view of the
absence of stirring. Moreover, we shall show that, up to a constant
the reverse inequality $ \l_0(N,\o) \leq C(1-\r)$ holds as well
for some choices of $\o$, see
the remark after the end of the proof. Similar results had been
obtained in \cite{Q1} in the homogeneous case $p_i\equiv\frac12$.

\begin{Th}
\la{main}
Assume (\ref{pi}) for some $\e>0$. For $m>1$, $\g=1$,
there exists $c_\e>0$ such that for any $N\geq 2$ and any $\o$ such that $0< \r <1$: 
\be
\l_1(N,\o)\geq c_\e\,.
\la{gap2}
\end{equation}
For $m>1$, $\g=0$, there exists $c_\e>0$ such that for any $N\geq 2$ and any $\o$ with $0<\r<1$:
\be
\l_0(N,\o)\geq c_\e\,(1-\r) \,.
\la{gap3}
\end{equation}
\end{Th}

\proof 
We start with some preliminary facts.
Consider functions $f$ of the occupation variables $\psi=\{\psi_i\}$ only.
Let $\cH_0$
denote the space of all such functions and observe that
$\cL_\g\cH_0\subset\cH_0$, i.e.\ $\cH_0$ is invariant, for both 
$\g=0,1$. This follows from
the fact that the only dependence of the rates on the configuration
$\eta$ is through the variables $\{\psi_i\}$, see (\ref{rates}). In particular, under the
generator
$\cL_\g$, the variables $\{\psi_i\}$ evolve as the Markov chain
of the case $m=1$. Therefore, the spectral gap estimate of Theorem
\ref{BEgg} applies to functions in $\cH_0$, for both $\g=0,1$.

For any $f\in L^2(\nu)$ we may write $f=f_0+f_0^\perp$, where $f_0=\nu(f\tc\cF_\psi)\in\cH_0$
and $f_0^\perp=f-f_0\in\cH_0^\perp$. Here $\cF_\psi$ denotes the
$\si$--algebra generated by the functions $\psi_i$ and 
$\cH_0^\perp$ is the orthogonal
complement of $\cH_0$, i.e.\ the space of $f$ such that
$\nu(f g)=0$ 
for all $g\in\cH_0$.
Since $\cL_\g$ leaves $\cH_0$ invariant, by self--adjointness it follows
that
$\cL_\g\cH_0^\perp\subset\cH_0^\perp$. In conclusion 
$$
\nu(f\cL_\g f) = \nu(f_0\cL_\g f_0) +\nu(f_0^\perp\cL_\g f_0^\perp) \,.
$$
Moreover, $\var_\nu(f) = \var_\nu(f_0) + \var_\nu(f_0^\perp)$. 
From these simple observations it is clear that 
the constant $\l_\g(N,\o)$
satisfies 
\be
\l_\g(N,\o)\geq \min\{\l^0_\g(N,\o),\l_{\g}^\perp(N,\o)\}\,,
\la{gammas}
\end{equation} 
where 
$\l_\g^0(N,\o),\l_\g^\perp(N,\o)$ are the 
constants obtained in the variational principle 
(\ref{gapkac}) -- applied to $\cL_\g$ --  by restricting
to $f\in\cH_0$ and $f\in\cH_0^\perp$ respectively.
From Theorem \ref{BEgg} we know that 
\be\la{goodda}
\l^0_\g(N,\o)\geq c_\e\,,
\end{equation}
for some $c_\e$, for both $\g= 0,1$. 
To prove Theorem \ref{main} we thus have to
estimate from below the constant $\l_\g^\perp(N,\o)$.

\subsubsection{The case $\g=1$} 
If $f\in\cH_0^\perp$, then
we must have $\nu(f\tc\cF_\psi)=0$.
Therefore
\be
\var_\nu(f)=\nu[\var(f\tc \cF_\psi)] + \var_\nu\left(\nu(f\tc \cF_\psi)\right)  = \nu[\var(f\tc \cF_\psi)] \,.
\la{varxi}
\end{equation}
Here $\var(\cdot\tc \cF_\psi)$ denotes variance w.r.t.\
$\nu(\cdot\tc\cF_\psi)$, 
the conditional probability given the 
color--blind configuration $\psi$, and we have used the standard decomposition of variance under conditioning. Let
$S=S(\eta):=\{i\in\{1,\dots,N\}:\;\,\psi_i=1\}$
denote the set of occupied sites. Clearly $|S|=\r N$, the total number of particles.  
Observe that once $\psi$ is known then the distribution 
of the variables $\eta$ is given by $\eta_i=0$
for $i\notin S$ (deterministically) 
and $\eta_i\in\{1,\dots,m\}$ on $S$ uniformly with
the constraint that $\sum_{i\in S}1_{\{\eta_i=k\}} = \o_k$,
$k\in \{1,\dots,m\}$. 
In particular there is no inhomogeneity once the set $S$ (or,
equivalently the configuration $\psi$) is
given. 
Therefore $\var(f\tc\cF_\psi )$ can be estimated for every
$\eta$
with the well known
bound
for random transpositions without disorder (see e.g.\ \cite{C1}):
\be
\var(f\tc \cF_\psi)\leq 
\frac1{4\,|S|}\sum_{i,j\in S}
\nu\left[(\grad_{ij}f)^2\tc \cF_\psi\right]
\,.
\la{varbo}
\end{equation} 
Here we use the notation $\grad_{ij}f(\eta) = f(\eta^{b})-f(\eta)$, $b=\{i,j\}$,
for the exchange gradient. 
Averaging with $\nu$ and using (\ref{varxi}) we obtain
(since $|S|=\r\, N$, deterministically)
\be
\var_\nu(f)\leq \frac1{4\,\r\,N}\,\,\,\nu\left(
\sum_{i,j\in S}(\grad_{ij}f)^2\right)\,.
\la{varbo2}
\end{equation} 
Suppose first $\r\geq \frac12$. 
Then (\ref{varbo2}), (\ref{diri}) and a uniform
lower bound on the rates $c_{b}$ imply
$$
\var_\nu(f)\leq \frac1{2\,N}\,\,
\sum_{i,j}\nu\left[(\grad_{ij}f)^2\right]
\leq C_\e 
\,\nu(f(-\cL_1)f)\,,
$$
for some constant $C_\e<\infty$, with the sum now extending to all pairs $i,j$. 
This shows that, with $c_\e=1/C_\e$: 
\be
\l_1^\perp(N,\o) \geq c_\e\,
,\quad\;\r\geq \frac12\,.
\la{lperp}
\end{equation}
Note that here we have used
$\g=1$ (if $\g=0$ there is no uniform lower bound on the rates $c_b$).

\smallskip

We turn to the proof in the case $\r\leq \frac12$. 
We rewrite (\ref{varbo2}) as
\begin{align}
\var_\nu(f) &\leq \frac1{4\,\r\,N}\sum_{i,j}
\nu\left[(\grad_{ij}f)^2
\,
1_{\{i\in S\}}1_{\{j\in S\}}
\right]\nonumber\\
& = \frac1{4\,\r(1-\r)N^2}
\sum_{i,j,\ell}
\nu\left[(\grad_{ij}f)^2
\,
1_{\{i\in S\}}1_{\{j\in S\}}1_{\{\ell\notin S\}}
\right]
\,,\la{os}
\end{align}
where we use the identity $(1-\r)N = N - \sum_{k=1}^m\o_k = 
\sum_{\ell}1_{\{\ell\notin S\}}$ for the number of empty sites.
Let $\eta\in\O_{N,\o}$ be fixed. 
Pick $i,j\in S(\eta)$ and write 
$\eta^{i,j}$ for the exchanged
configuration.
Observe that for any $\ell\notin S(\eta)$ we
can write
$$
\eta^{i,j} = \left(\left(\eta^{i,\ell}\right)^{i,j}\right)^{j,\ell}\,.
$$
Therefore \begin{align*}
\grad_{ij}f(\eta) &= [f(\eta^{i,j})-f(\eta)]\\
&= \grad_{j\ell}f\left(\left(\eta^{i,\ell}\right)^{i,j}\right)
+ \grad_{ij}f\left(\eta^{i,\ell}\right)
+\grad_{i\ell}f\left(\eta\right)
\,.
\end{align*}
Each term in the sum appearing in (\ref{os}) can then be estimated by
\begin{align}
 &\nu\left[(\grad_{ij}f)^2\,1_{\{i\in S\}}1_{\{j\in S\}}1_{\{\ell\notin S\}}
\right]
\la{os0}\\&\leq 3\,\nu
\left[\left\{\left(\grad_{j\ell}f
\left(\left(\eta^{i,\ell}\right)^{i,j}\right)\right)^2\,
+\left(\grad_{ij}f\left(\eta^{i,\ell}\right)\right)^2
+\left(\grad_{i\ell}f\right)^2\right\}
1_{\{i\in S\}}1_{\{j\in S\}}1_{\{\ell\notin
  S\}}\right] \,.
\nonumber
\end{align}
Next, we claim that there exists $C_1 = C_1(\e)<\infty$ such that
\begin{align}
&\nu
\left[\left(\grad_{j\ell}f
\left(\left(\eta^{i,\ell}\right)^{i,j}\right)\right)^2
1_{\{i\in S\}}1_{\{j\in S\}}1_{\{\ell\notin
  S\}}\right] \nonumber\\
&\quad\quad\quad\qquad\;\;\leq C_1 \,\nu
\left[\left(\grad_{j\ell}f(\eta)\right)^2
1_{\{i\in S\}}1_{\{j\notin S\}}1_{\{\ell\in
  S\}}\right]\,.
\la{os1}
\end{align}
Note the change of the indicator functions in (\ref{os1}).
To prove (\ref{os1}) we write, with the change of variables
$\varphi:=\eta^{i,\ell}$ and $\b:=\varphi^{i,j}$
\begin{align*}
&\nu
\left[\left(\grad_{j\ell}f
\left(\left(\eta^{i,\ell}\right)^{i,j}\right)\right)^2
1_{\{i\in S\}}1_{\{j\in S\}}1_{\{\ell\notin
  S\}}\right]\\
&\qquad\qquad = \sum_{\eta}\nu(\eta)\,\left(\grad_{j\ell}f
\left(\left(\eta^{i,\ell}\right)^{i,j}\right)\right)^2\,
1_{\{i\in S(\eta)\}}1_{\{j\in S(\eta)\}}1_{\{\ell\notin
  S(\eta)\}}\\
&\qquad\qquad =\sum_{\varphi}\nu(\varphi^{i,\ell})\,\left(\grad_{j\ell}f
\left(\varphi^{i,j}\right)\right)^2\,
1_{\{i\notin S(\varphi)\}}1_{\{j\in S(\varphi)\}}1_{\{\ell\in
  S(\varphi)\}}\\
&\qquad\qquad =\sum_{\b}\nu((\b^{i,j})^{i,\ell})\,
\left(\grad_{j\ell}f
\left(\b\right)\right)^2\,
1_{\{i\in S(\b)\}}1_{\{j\notin S(\b)\}}1_{\{\ell\in
  S(\b)\}}\\
&\qquad\qquad= \nu
\left[\chi_{i,j,\ell}\,\left(\grad_{j\ell}f\right)^2
1_{\{i\in S\}}1_{\{j\notin S\}}1_{\{\ell\in
  S\}}\right]\,,
\end{align*}
where $$
\chi_{i,j,\ell}(\eta):=\frac{\nu((\eta^{i,j})^{i,\ell})}{\nu(\eta)}\,,
$$
is the change of measure. 
Since the variables $p_i$ defining our measure satisfy (\ref{pi})
it is not hard to check that $\chi_{i,j,\ell}(\eta)\leq C_\e$ uniformly for
some constant $C_\e$. This proves (\ref{os1}).

\smallskip
Moreover, in a similar way one proves that there is a constant
$C_2=C_2(\e)<\infty$ such that  
\begin{align}
&\nu
\left[\left(\grad_{ij}f
\left(\eta^{i,\ell}\right)\right)^2
1_{\{i\in S\}}1_{\{j\in S\}}1_{\{\ell\notin
  S\}}\right] \nonumber\\
&\quad\quad\quad\qquad\;\;\leq C_2 \,\nu
\left[\left(\grad_{ij}f(\eta)\right)^2
1_{\{i\notin S\}}1_{\{j\in S\}}1_{\{\ell\in
  S\}}\right]\,.
\la{os2}
\end{align}
From (\ref{os}) and (\ref{os0}), 
using (\ref{os1}) and (\ref{os2}) we obtain 
for a suitable constant $C_3$:
\begin{align}
&\var_\nu(f) \leq \frac{C_3}{\r(1-\r)N^2}
\sum_{i,j,\ell}
\Big\{\nu\left[(\grad_{j\ell}f)^2
\,
1_{\{i\in S\}}1_{\{j\notin S\}}1_{\{\ell\in S\}}
\right]\nonumber\\
&+ \nu\left[(\grad_{ij}f)^2
\,
1_{\{i\notin S\}}1_{\{j\in S\}}1_{\{\ell\in S\}}
\right] + 
\nu\left[(\grad_{i\ell}f)^2
\,
1_{\{i\in S\}}1_{\{j\in S\}}1_{\{\ell\notin S\}}
\right]\Big\}\,.
\la{os3}
\end{align}
Since $\sum_{i}1_{\{i\in S\}}=\r N$, setting $C_4=3\,C_3$ 
we see that (\ref{os3}) can be written as
\be
\var_\nu(f) \leq \frac{C_4}{(1-\r)N}
\sum_{i,j}
\nu\left[(\grad_{ij}f)^2
\,
1_{\{i\in S\}}1_{\{j\notin S\}}
\right]\,.
\la{good}
\end{equation}
Using $\r\leq \frac{1}2$ and (\ref{diri}) we see that
$$
\var_\nu(f) \leq C\,\nu(f(-\cL_1)f)\,,
$$
for some constant $C=C(\e)<\infty$ and for all $f\in\cH_0^\perp$. 
Therefore 
we have proved that $\l_1^\perp(N,\o)\geq c_\e\,$. 
Together with
(\ref{gammas}), (\ref{goodda}) and (\ref{lperp}) 
we see that $\l_1(N,\r)$ is uniformly
bounded from below and the claim (\ref{gap2}) holds. This proves Theorem
\ref{main} in the case $\g=1$.

\subsection{The case $\g=0$}
Here we cannot use the argument leading to (\ref{lperp})
above. However,
we can use the argument giving (\ref{good}) without modification
(it never used the fact that $\g=1$). 
In particular, (\ref{good}) holds 
for any $0<\r <1$ and any $f \in\cH_0^\perp$.
%
Now, observe that for any edge $b=\{i,j\}$ such that 
$i\in S$ and $j\notin S$ the rate $c_b^0$ defined in (\ref{rates})
is uniformly bounded away from
zero with constants only depending on the $\e$ from
(\ref{pi}) (this follows from the fact that for such cases either
$\psi_i(1-\psi_j)=1$ or $\psi_j(1-\psi_i)=1$). Therefore, there exists $C=C(\e)<\infty$ such that
for any $f \in\cH_0^\perp$:
\be
\var_\nu(f) \leq \frac{C}{(1-\r)}
\,\nu(f(-\cL_0)f)\,.
\la{good2}
\end{equation}
This proves that $\l_0^\perp(N,\o)\geq c\,(1-\r)$. 
From
(\ref{gammas}) and (\ref{goodda}) we see that 
$\l_0(N,\o)$ satisfies the claim (\ref{gap3}). This proves Theorem
\ref{main} in the case $\g=0$. \qed

\bigskip

Let us briefly address the issue of 
comparable upper bounds on spectral gaps. 
For example, to prove that  
$\l_0(N,\o)=O(1-\r) $, as $\r\to 1$ 
one can 
consider the case of $m=2$ colors, with $\o_1=\o_2$ 
such that $\o_1+\o_2 = \r\,N$
with $\r\to 1$. Then choose
$f$ in the variational principle
(\ref{gapkac}) as the indicator function of the event that 
a given vertex $i$ is occupied by a particle of color $1$. 
The variance of $f$ is of order $1$. On the other hand, since the
number of empty sites that can be used to change the value of $\eta_i$ is
$(1-\r)N$,
it is not hard to show that 
$\nu(f(-\cL_0)f)$ is of order $(1-\r)$. It follows that
$\l_0(N,\o)=O(1-\r)$. Of
course, if $\g=1$ then $\nu(f(-\cL_1)f)$ is of order $1$ and the gap
does not vanish as $\r\to 1$ in that case. 

\smallskip

Finally, we point out an interesting problem concerning 
local versions of the exclusion 
dynamics described by (\ref{diri}). The local dynamics is obtained 
by summing over pairs $b$ which are given by the edges of a small
subgraph of the complete graph,
such as e.g.\ the box of side $L\sim N^{1/d}$ in a $d$--dimensional grid: 
$\L_L=[1,L]^d\cap\bbZ^d$. In the latter cases one expects the gap to
scale as $L^{-2}$. In the case $m=1$ there is a nice argument 
in \cite[Lemma 5.1 and Lemma 5.2]{Q} which allows to derive
such an estimate from the complete--graph bound (\ref{exc1}). On
the other hand, the case $m>1$ with $\g=0$ 
seems to be much more complicated and
we are not aware of any result in this direction, except for the
homogeneous case considered in \cite{Q1}.

%
%
%


 \smallskip
\bigskip

\noindent
{\bf Acknowledgments.}
I would like to thank Eric A.\ Carlen, Maria C.\ 
Carvalho, Paolo Dai Pra, Gustavo Posta 
and Prasad Tetali for several useful discussions around this work.

\end{document}